\documentclass[centertags,leqno]{article}

\usepackage{amssymb}
\usepackage{latexsym}
\usepackage{amsxtra}
\usepackage{amscd}
\usepackage{theorem}
\usepackage{xypic}

\theoremstyle{plain}

{\theorembodyfont{\slshape}

        \newtheorem{thm}{Theorem}[section]
        \newtheorem{cor}[thm]{Corollary}
        \newtheorem{lem}[thm]{Lemma}
        \newtheorem{prop}[thm]{Proposition}

}

{\theorembodyfont{\rmfamily}

        \newtheorem{defn}[thm]{Definition}
        
        \newtheorem{rem}[thm]{Remark}
        \newtheorem{exa}[thm]{Example}

}

\renewcommand{\em}{\sl}

\newcommand{\proof}{{\bf Proof:\ }}
\newcommand{\Endproof}{\hspace*{\fill} $\Box$ \vspace{1ex} \noindent }

\makeatletter
\renewcommand{\subsection}{\@startsection{subsection}{2}%
        {\z@}{-3.25ex plus -1ex minus-.2ex}{-1em}{\bf}}
\makeatother

\renewcommand{\arraystretch}{1.5}

\newcommand{\PP}{\mathbb{P}}
\newcommand{\ZZ}{\mathbb{Z}}

\newcommand{\QQ}{\mathbb{Q}}

\newcommand{\Nc}{{\mathcal N}}
\newcommand{\FF}{\mathbb{F}}
\renewcommand{\AA}{\mathbb{A}}

\newcommand{\T}{\mathcal{T}}
\newcommand{\E}{\mathcal{E}}
\newcommand{\M}{\mathcal{M}}
\newcommand{\LL}{\mathcal{L}}
\newcommand{\OO}{\mathcal{O}}

\newcommand{\olog}{\Omega^{\rm\scriptscriptstyle log}}
\newcommand{\tlog}{\tau^{\rm\scriptscriptstyle log}}

\newcommand{\diff}{{\rm d}}
\newcommand{\Fil}{{\rm F}^1}
\newcommand{\Gr}{{\rm Gr}^1}

\newcommand{\Hom}{\mathop{\rm Hom}\nolimits}
\newcommand{\End}{\mathop{\rm End}\nolimits}

\newcommand{\GL}{\mathop{\rm GL}}
\newcommand{\m}{\mathfrak{m}}
\newcommand{\ord}{{\rm ord}}
\newcommand{\ords}{{\rm \scriptscriptstyle ord}}

\newcommand{\PSL}{\mathop{\rm PSL}\nolimits}

\newcommand{\inj}{\hookrightarrow}
\newcommand{\To}{\longrightarrow}
\newcommand{\iso}{\stackrel{\sim}{\to}}

\newcommand{\lpfeil}[1]{\stackrel{#1}{\To}}

\newcommand{\Spec}{\mathop{\rm Spec}\nolimits}

\newcommand{\bsigma}{{\boldsymbol\sigma}}

\title{Indigenous bundles with nilpotent $p$-curvature}
\author{Irene I.\ Bouw and Stefan Wewers}
\date{}

\begin{document}
\maketitle

\begin{abstract}
  We study indigenous bundles in characteristic $p>0$ with nilpotent
  $p$-curvature, and show that they correspond to so-called
  deformation data. Using this equivalence, we translate the existence
  problem for deformation data into the existence of polynomial
  solutions of certain differential equations with additional
  properties. As in application, we show that $\PP^1$ minus four
  points is hyperbolically ordinary (in the sense of Mochizuki
  \cite{Mochizuki1}).  We also give a concrete application to
  existence of deformation data with fixed local invariants.

\vspace{1ex}
\noindent {\sl\small Mathematics Subject Classification} (2000): Primary
14H25, 14H60.
\end{abstract}

\section{Introduction}\label{introsec}
The main goal of this paper is to set up an equivalence between
indigenous bundles with nilpotent $p$-curvature and deformation data.
Deformation data arise in the theory of stable reduction of Galois
covers of curves. An aspect of this theory which has not yet been
treated in a satisfactory way is the existence of deformation data
with given local invariants. Though deformation data are in principal
easier objects than indigenous bundles, the existence problem seems
more approachable via indigenous bundles. This is especially true for
indigenous bundles on $\PP^1$.  Here we give a concrete interpretation
of an indigenous bundle in terms of a solution of an ordinary
differential equation in characteristic $p$ with additional
properties. The existence problem of deformation data becomes then a
variant of Dwork's accessory parameter problem (\cite{Dwork}).

Let $X$ be a smooth projective curve over an algebraically closed field $k$ of
characteristic $p>0$.  A deformation datum $(Z, \omega)$ consist of a cyclic
cover $Z\to X$ of order dividing $p-1$ together with a 
differential form $\omega$ on $Z$, satisfying certain additional properties
(Definition \ref{defodatdef}).  Associated to a
deformation datum is a set of local invariants called the signature.

Let $f_K:Y_K\to X_K$ be a $G$-Galois cover defined over a field $K$ of
characteristic zero. Suppose that $p$ strictly divides the order of
$G$, and that $f_K$ has bad reduction at some place $\wp$ of $K$ with
residue characteristic $p$. Let $\bar{f}$ denote the {\em stable
reduction} of ${f}_K$ at $\wp$ (see \cite{bad} for a precise defintion). Then
to $\bar{f}$ one may associate a set of deformation data; they
describe the inseparable part $\bar{f}$ of $f_K$. A particularly
interesting case is when $f_K$ is a Belyi map, i.e.\ when
$X_K=\PP^1_K$ and $f_K$ is branched at $\{0,1,\infty\}$.  Here one has
a good description of the stable reduction $\bar{f}$, see \cite{bad},
\cite{Raynaud98}. In particular, the inseparable part of $\bar{f}$ is
described by one deformation datum on $X:=\PP^1_k$, where
$k=\bar{\FF}_p$. Deformation data coming from the stable reduction of
Belyi maps are called {\sl special}; they are characterized by a
numerical condition on the signature.

The correspondence between Belyi maps with bad reduction and special
deformation data is intriguing and, at present, quite mysterious. The
only case that is fully understood is worked out in
\cite{mcav}. Suppose $f_K:Y_X\to X_K=\PP^1_K$ is a Belyi map which
is Galois and whose Galois group is contained in the linear group
${\rm GL}_2(\FF_p)$. If $f_K$ has bad reduction at $\wp$, then the
corresponding special deformation datum is {\em hypergeometric}, i.e.\
it corresponds to a polynomial solution of a certain hypergeometric
differential equation. Furthermore, the correspondence between the
group theoretic description of the cover $f_K$ and the parameters of
this hypergeometric equation is totally explicit. As an
application, one gets an if and only if criterion for a Belyi map with
Galois group contained in ${\GL}_2(\FF_p)$ to have good reduction at
$p$. 

Deformation data show up naturally in other situations, too. One
example is the {\em lifting problem}, i.e.\ the question of whether a
curve with given automorphism group in characteristic $p$ lifts to
characteristic $0$. See e.g.\ \cite{Dp}. Here, as in other situations
where deformation data play a role, the difficult problem is to show
the existence of deformation data with a given signature. One of the
motivations behind the present paper is to translate this existence
problem into the existence of a different kind of object which is
often more tractable.

In \S \ref{defodat} we associate to every deformation datum
$(Z,\omega)$ over $X$ the equivalence class of a certain flat vector
bundle of rank two $(\E,\nabla)$ on $X$, with at most regular
singularities. (The equivalence class corresponds to the associated
projective bundle $\PP(\E)$, together with its induced connection.)
We show that $(\E,\nabla)$ is an {\em indigenous bundle} in the sense
of Mochizuki, see \cite{Mochizuki1}, \cite{Mochizuki2}. Moreover,
$(\E,\nabla)$ has nilpotent and nonvanishing
$p$-curvature. Conversely, we show that every indigenous bundle with
nonvanishing and nilpotent $p$-curvature comes from a deformation
datum $(Z,\omega)$. A technical part is to relate the invariants on
both sides of the equivalence. We show that the signature of a
deformation datum can be expressed in terms of the order of the zeros
of the $p$-curvature of the corresponding indigenous bundle.  The
correspondence between indigenous bundles with nilpotent $p$-curvature
and deformation data is implicit in Mochizuki's work.  Our
construction is more direct and explicit, and is inspired by work of
Ihara \cite{Ihara74}.

The last part of the paper concerns the existence of deformation data
on $X=\PP^1_k$ with given signature. In this case, there is a
classical correspondence between indigenous bundles and second order
differential equations
\[
   L(u)=u''+p_1u'+p_2u=0
\] 
with at most regular singularities. Indigenous bundles with nilpotent
$p$-curvature correspond to differential equations which admit a
polynomial solution. See \S \ref{accsec} for a more precise
statement. Using the results of the first part of the paper, we
translate the existence problem for deformation data into a variant of
Dwork's accessary parameter problem (\cite{Dwork}).  We use this
description to discuss the existence of special deformation data with
four singularities in \S \ref{specialsec}.

As an application of the description of indigenous bundles on $\PP^1$
in terms of solutions of a differential equation, we show in \S
\ref{ordinarysec} that a projective line minus four points is
hyperbolically ordinary (\cite{Mochizuki1}). This complements a result
of Mochizuki which states that every generic marked curve is
hyperbolically ordinary. In terms of deformation data our result
implies that for every $\lambda\in k-\{0,1\}$ there exists a
deformation datum $(Z,\omega)$ with singularities
$\{0,1,\infty,\lambda\}$ and signature
${\boldsymbol\sigma}=(0,0,0,0)$.

We also briefly discuss the analogous question for other signatures
${\boldsymbol\sigma}$. In \S \ref{ordtorsec} we give a numerical condition on
${\boldsymbol\sigma}$ which determines the dimension of the space of
deformation data on $\PP^1$ with four singular points and signature
${\boldsymbol\sigma}$. It is either zero or one, and both cases
occur. In particular, for special deformation data it is zero.  It
therefore seems hard to prove the existence of special deformation
data by using Mochizuki's theory, which is best suited to prove
existence results on generic curves. A different approach to
construct special deformation data, which also uses the language of
indigenous bundles, is developed in \cite{IreneHabil}.


\section{Indigenous bundles}

\subsection{}

The following notation will be fixed throughout this paper.  Let $k$
be an algebraically closed field and $X$ a smooth projective connected
curve of genus $g$ over $k$. We fix $r\geq 0$ pairwise distinct closed
points $x_1,\ldots,x_r\in X$, which we call {\em marked points}. We
assume that $2g-2+r>0$. We denote by $\olog_{X/k}=\Omega_{X/k}(\sum
x_i)$ the sheaf of differential $1$-forms on $X$ with at most simple
poles in the marked points and $\tlog_{X/k}\cong(\olog_{X/k})^{-1}$
its dual, i.e.\ the sheaf of vector fields on $X$ with at least simple
zeros in the marked points.

A {\em flat vector bundle} is a vector bundle $\E$ on $X$ together
with a connection $\nabla:\E\to\E\otimes\olog_{X/k}$. This means that
the connection has regular singularities at the marked points (see
e.g.\ \cite{Katz}, \S 11). Two flat vector bundles $\E$ and $\E'$ are
called {\em equivalent} if there exists a flat line bundle $\LL$ and a
horizontal isomorphism $\E'\cong\E\otimes\LL$. Notation: $\E\sim\E'$.

\subsection{}\label{monosec}

Let $(\E,\nabla)$ be a flat vector bundle on $X$ of rank two. For $i=1,..,r$,
we define the {\em monodromy operator} $\mu_i$ as an endomorphism of the fiber
$\E|_{x_i}$ of $\E$ at $x_i$, as follows. Let $t$ be a local parameter at
$x_i$. Then $\nabla(t\partial/\partial t)$ defines a $k$-linear endomorphism
of the stalk $\E_{x_i}$ of $\E$ at $x_i$ which fixes the submodule
$\m_{x_i}\cdot\E_{x_i}$. Here $\m_{x_i}$ denotes the maximal ideal of the
local ring $\OO_{X,x_i}$. Therefore, $\nabla(t\partial/\partial t)$ induces a
$k$-linear endomorphism $\mu_i$ of the fiber
$\E|_{x_i}=\E_{x_i}/\m_{x_i\cdot}\E_{x_i}$. One checks easily that $\mu_i$
does not depend on the choice of the parameter $t$.

Let $\alpha_i,\beta_i$ be the two eigenvalues of $\mu_i$. We call
$\alpha_i,\beta_i$ the {\em local exponents} of $\nabla$ at $x_i$. We
distinguish two cases. If $\mu_i$ is not semisimple, then
$\alpha_i=\beta_i$ and 
\[
   \mu_i \sim 
     \begin{pmatrix} \alpha_i & 1 \\ 0 & \alpha_i \end{pmatrix}.
\]
If this is the case then we say that $\nabla$ has {\em logarithmic
monodromy} at $x_i$. If $\mu_i$ is semisimple then 
\[
   \mu_i \sim 
     \begin{pmatrix} \alpha_i & 0 \\ 0 & \beta_i \end{pmatrix},
\]  
and we say that $\nabla$ has {\em toric monodromy} at $x_i$.

\subsection{}

A {\em filtration} on a flat vector bundle $(\E,\nabla)$ of rank two
consists of a line subbundle $\Fil\E\subset\E$ such that
$\Gr\E:=\E/\Fil\E$ is also a line bundle. For such a filtration, the
connection $\nabla$ induces a {\em Kodaira--Spencer} map
\[
      \kappa:\Fil\E \To \Gr\E\otimes\olog_{X/k}.
\]
If it seems more convenient, we  regard $\kappa$
as a morphism
\[
      \kappa:\tlog_{X/k} \To (\Fil\E)^{-1}\otimes\Gr\E.
\]
Note that, written in either way, $\kappa$ is $\OO_X$-linear. 

\begin{defn}\label{indidef}
  An {\em indigenous bundle} on $X$ is a flat vector bundle
  $(\E,\nabla)$ of rank two which satisfies the following conditions:
  \begin{enumerate}
  \item
    There exists a filtration $\Fil\E\subset\E$  whose
    associated Kodaira--Spencer map is an isomorphism.
  \item
    The connection $\nabla$ has nontrivial monodromy at every marked
    point (i.e.\ $\mu_i\not=0$ for all $i$).
  \end{enumerate}
  A filtration $\Fil\E\subset\E$ as in (i) is called a {\em Hodge
  filtration}.
\end{defn}

\begin{prop} \label{indiprop1}
  Let $(\E,\nabla)$ be an indigenous bundle. 
  \begin{enumerate}
  \item
    The Hodge filtration for $(\E,\nabla)$ is unique.
  \item 
    Fix an index $i\in\{1,\ldots,r\}$ and suppose that $\nabla$ has toric
    monodromy at $x_i$, with exponents $\alpha_i$, $\beta_i$. Then
    $\alpha_i\not=\beta_i$. 
  \item
    Any equivalent bundle $\E'\sim\E$ is indigenous as well. 
 \end{enumerate}
\end{prop}

 Note that Part (i) of the proposition implies that it makes sense to
  speak about {\em the} Hodge filtration of an indigenous bundle.

\bigskip
\proof Define the {\em height} of a filtration $\Fil\E\subset\E$ as
the integer
\[
       \frac{1}{2}(\deg\Gr\E - \deg\Fil\E).
\]
By assumption, a Hodge filtration has height $-g+1-r/2<0$. On the other hand,
it is easy to see that two filtrations with negative height are equal. This
proves (i). To prove (ii), suppose that $\alpha_i=\beta_i$. Then $\mu_i$ fixes
the line $\Fil\E|_{x_i}\subset\E|_{x_i}$. But this would mean that the
Kodaira--Spencer map $\kappa$ vanishes at $x_i$, which gives a contradiction.
The proof of (iii) is easy and left to the reader.  \Endproof

\begin{rem}
  In \cite{Mochizuki1}, indigenous bundles are required to have
  logarithmic monodromy at the marked points and to have trivial
  determinant. What we call here an indigenous bundle is a flat vector
  bundle whose associated projective bundle is  {\em
  torally indigenous},  in the sense of \cite{Mochizuki2}, \S I.4.
\end{rem}

\section{Indigenous bundles in characteristic $p$}

\subsection{}\label{indidefsec}

From now on, we assume that our base field $k$ has odd, positive
characteristic $p$. In this section, we develop the theory of
indigenous bundles in characteristic $p>2$. All definitions are
essentially due to Mochizuki (\cite{Mochizuki1}, \cite{Mochizuki2}).

We set $\T:=(\tlog_{X/k})^{\otimes p}$. This is a
line bundle on $X$ of degree $-p(2g-2+r)<0$. We endow $\T$ with the
unique connection $\nabla_{\T}:\T\to\T\otimes\olog_{X/k}$ such that
the subsheaf $\T^\nabla$ of horizontal sections consists precisely of
the `$p$-th powers', i.e.\ of sections of the form $D^{\otimes p}$,
where $D$ is a section of $\tlog_{X/k}$ (\cite{Katz}, Theorem 5.1).

Let $(\E,\nabla)$ be a flat vector bundle on $X$. The {\em
$p$-curvature} of $(\E,\nabla)$ is an $\OO_X$-linear morphism
\[
     \Psi_\E:\T \To \End_{\OO_X}(\E),
\]
defined as follows. Let $D$ be a rational section of $\tlog_{X/k}$. We
regard $D$ as a derivation of the function field $K=k(X)$. Then
$D^p:=D\circ\cdots\circ D$ is again a derivation of $K$ and
$\nabla(D)$ and $\nabla(D^p)$ are $k$-linear endomorphisms of the
$K$-vector space $E:=\E\otimes_{\OO_X}K$. We define
\[
        \Psi_\E(D^{\otimes p}) := \nabla(D)^p - \nabla(D^p).
\]
This is a $K$-linear endomorphism of $E$.  One shows that the rule
$D^{\otimes p}\mapsto\Psi(D^{\otimes p})$ descents to
the desired $\OO_X$-linear map $\Psi_\E$ (\cite{Katz}, 5.0.1). 

It is important to notice that the $p$-curvature is {\em horizontal}
in the sense that it commutes with the canonical connections on $\T$
and $\End_{\OO_X}(\E)$. Indeed, by the definitions of these connections,
the claim that $\Psi_\E$ is horizontal is equivalent to the fact that the
endomorphisms $\Psi_\E(D^{\otimes p})$ and $\nabla(D)$ of $E$
commute. This is easy to check, see also \cite{Katz}, 5.2.2.    

\begin{defn} \label{nilpdef}
  An indigenous bundle $(\E,\nabla)$ on $X$ is called
  \begin{enumerate}
  \item
    {\em active} if $\Psi_\E\not=0$, 
  \item
    {\em nilpotent} if the image of $\Psi_\E$ consists of nilpotent
    endomorphisms.
  \end{enumerate}
  
\end{defn}

Let us, from now on, assume that $(\E,\nabla)$ is active and
nilpotent. We are  only interested in
the equivalence class of $\E$ (i.e.\ in the associated projective
bundle $\PP(\E)$). Therefore we may replace $\E$ by any equivalent bundle $\E'$
which is also active and nilpotent. 

\begin{defn}
  We say that $(\E,\nabla)$ is {\em normalized} if there exists a
  horizontal and surjective homomorphism $\gamma:\E\to\OO_X$. 
\end{defn}

Normalized indigenous bundles correspond to Mochizuki's ${\rm FL}$-bundles
(\cite{Mochizuki1}, \S 2.1).

We claim that for any indigenous bundle $(\E,\nabla)$ which is active
and nilpotent there is an equivalent bundle
$(\E',\nabla')$ which is normalized. Moreover,
$(\E',\nabla')$ is unique, up to isomorphism. 

Let $\M\subset\E$ be the kernel of $\Psi_\E$, i.e.\ the
maximal subbundle on which $\Psi_\E$ is zero. Our assumption implies
that $\M$ and the quotient $\LL:=\E/\M$ are line bundles. Moreover, it
follows from the fact that $\Psi_\E$ is horizontal that $\M$ is
invariant under the connection $\nabla$. In other words, we obtain a
short exact sequence of flat vector bundles
\begin{equation}\label{exactseqeq}
      0 \;\to\; \M \To \E \To \LL \;\to\; 0.
\end{equation}
The $p$-curvature of the induced connections on $\M$ and
$\LL$ is zero.

We set $\E':=\E\otimes\LL^{-1}$. Since the $p$-curvature on $\LL$ is
zero, $\LL^{-1}$ admits a nonzero rational horizontal section. Using
this fact, it is easy to see that the $p$-curvature of $\E'$ is equal
to the $p$-curvature of $\E$ (here we identify $\End(\E)$ and
$\End(\E')$ in the obvious way). In particular, $\E'$ is active and
nilpotent. By construction, $\OO_X$ is a quotient of $\E'$. Hence
$\E'$ is normalized. The uniqueness of $\E'$ follows from the fact
that $\M$ is the maximal saturated line subbundle of $\E$ which is
invariant under $\nabla$.

\subsection{} \label{pcurvsec}

From now on, we  assume that $\E$ is normalized.  Moreover, we
 fix a surjective and horizontal morphism $\gamma:\E\to\OO_X$.
Note that $\gamma$ is unique up to  multiplication by an element of
$k^\times$. By definition,  $\M$ is the kernel of $\gamma$.

\begin{defn}\label{spikedef}
A point $x\in X$ where $\Psi_\E$ vanishes is
  called a {\em spike}. We write $n_x:=\ord_x(\Psi_\E)$ for the order of
  vanishing of $\Psi_\E$ at $x$ and say that $x$ is a {\em spike of order
  $n_x$}.
\end{defn}

 We may regard
the $p$-curvature of $\E$ as a nonzero horizontal homomorphism
\[
   \Psi_\E:\T\To \Hom_{\OO_X}(\LL, \M)=\M\otimes \LL^{-1}\simeq \M,
\]
since we assumed that $\M$ is normalized.
Let $S=\sum_x n_x\cdot x$ be the divisor of zeros of $\Psi_\E$ (with
support in the spikes). Then $\Psi_\E$ induces an isomorphism
$\T(S)\cong\M$. We let $\delta:\T(S)\inj\E$ denote the composition of
the last isomorphism {\em multiplied by $-1$} with the canonical
injection $\M\inj\E$. From now on, we identify $\T(S)$ with a
subbundle of $\E$ via the injection $\delta$. Since  $\E$ is normalized, the
exact sequence (\ref{exactseqeq}) becomes
\[
   0 \;\to\; \T(S) \To \E \;\lpfeil{\gamma}\; \OO_X \;\to\; 0.
\]
It follows that the $p$-curvature of $\E$, considered as $\OO_X$-linear
morphism $\T\to \End_{\OO_X}(\E)$, is given by the
formula
\begin{equation} \label{FLeq2}
   \Psi_\E(D^{\otimes p})(e) = -\gamma(e)\cdot D^{\otimes p}.
\end{equation} 
Here $e$ is a  section of $\E$ and $D$ a vector field.

Let $e_1$ be a rational section of the line subbundle $\Fil\E$. We set
$u:=\gamma(e_1)\in\OO_{X,x}$ and $e_0:=u^{-1}e_1$. By definition,
$e_0$ is a rational section of $\Fil\E$ with $\gamma(e_0)=1$.  Note
that these two properties characterize $e_0$ uniquely. In particular,
$e_0$ is independent of the choice of $e_1$. (However, $e_0$ does
depend on the choice of $\gamma$.) We call $e_0$ the {\em canonical
section} of $\E$.

\begin{lem} \label{sslem}
  The canonical section $e_0$ has at most simple poles and no
  zeros. The marked points $x_i$ are  regular points of $e_0$.
\end{lem}

\proof Let $x\in X$ be a closed point and $t$ a local parameter at
$x$. We set $D:=t\partial/\partial t$ or $D:=\partial/\partial t$,
depending on whether $x$ is marked or not. We may assume that the
section $e_1$ is a generator of $\Fil\E$ at $x$. Then $u=\gamma(e_1)$
is an element of the local ring $\OO_{X,x}$. We have to show that $u$
has at most  simple zeros, and is invertible if $x$ is a marked point.

Let $e_2:=\nabla(D)(e_1)$. Since the Kodaira--Spencer morphism is
nonzero at $x$, $(e_1,e_2)$ is an $\OO_{X,x}$-basis for the stalk of
$\E$ at $x$. Moreover, we have $\gamma(e_2)=D(u)$.  Let $e_3$ be any
element of $\E_x$ with $\gamma(e_3)=1$ and write $e_3=q_1e_1+q_2e_2$,
with $q_1,q_2\in\OO_{X,x}$. Applying $\gamma$ to to this equality we
get
\[
    q_1u + q_2D(u) = 1,
\]
and hence
\[
      \min(\,\ord_xu, \ord_xD(u)\,) = 0.
\]
Suppose that $x$ is not marked. Then $D=\partial/\partial t$ and
$\ord_xD(u)\geq \ord_xu-1$. We conclude that $\ord_x(u)\in\{0,1\}$. On the
other hand, if $x=x_i$ is a marked point, then $D=t\partial/\partial t$ and
$\ord_xD(u)\geq\ord_xu$. Therefore, $\ord_xu=0$.  \Endproof

\begin{defn}\label{ssdef}
  The points where the canonical section $e_0$ has a simple pole are
  called the {\em supersingular points} of $(\E,\nabla)$. 
\end{defn}

\subsection{}\label{monosec2}
Let $(\E, \nabla)$ be an active, nilpotent indigenous  bundle which we assume
to be normalized.
We now compute the $p$-curvature in the marked points. Let $t$ be a local
parameter at the marked point $x_i$, and set $D_i:=t\partial/\partial t$. It
is easy to see that $D_i^p=D_i$. Therefore,
\[
   \Psi_{x_i} := \Psi_\E(D_i^{\otimes p})|_{x_i} =
     \mu_i^p-\mu_i.
\]

Since the subbundle ${\rm Ker}(\gamma)\subset\E$ is invariant under $\nabla$,
the monodromy operator $\mu_i$ fixes the line ${\rm
  Ker}(\gamma)|_{x_i}\subset\E|_{x_i}$.  Therefore, ${\rm Ker}(\gamma)$ has
regular singularities at $x_i$.   It follows that the local exponents of $(\E,
\nabla)$ at $x_i$ are $(\alpha_i, 0)$, where $\alpha_i$ is the local exponent
of $({\rm Ker}(\gamma), \nabla)$.

If $\E$ has a logarithmic singularity at $x_i$ with local exponents
$(\alpha_i, \alpha_i)=(0,0)$ then
\[
   \Psi_{x_i} \;\sim\; 
\begin{pmatrix} 0&-1\\0&0\end{pmatrix}.
\]
In particular, $\Psi_{x_i}\not=0$, so $x_i$ is  not a spike.  

Similarly, if $\E$ has a toric singularity at $x_i$, with local exponents
$(\alpha_i, \beta_i)=(\alpha_i,0)$, then
\[
   \Psi_{x_i} \;\sim\; \begin{pmatrix} 
         \alpha_i^p-\alpha_i & 0 \\ 0 & 0
      \end{pmatrix}.
\]
Since $(\E, \nabla)$ is nilpotent, we have that $\Psi_{x_i}^2=0$, This implies
that $\alpha_i\in\FF_p^\times$. In particular, $x_i$ is a spike.

\begin{prop} \label{nilpotprop1}
  Suppose that $(\E,\nabla)$ is active, nilpotent, and normalized. Let $x\in
  X$ be a closed point.
  \begin{enumerate}
  \item If $x=x_i$ is a marked point and $\nabla$ has logarithmic monodromy at
    $x$, then $n_x=0$.
  \item If $x=x_i$ is a marked point and $\nabla$ has toric monodromy at $x$,
    with exponents $(\alpha_i,0)$, then $x$ is a spike. The order $n_x$ of
    this spike satisfies
    \[
       n_x \;\equiv\;  - \alpha_i \;\not\equiv\; 0 \pmod{p}.
    \]
  \item
    If $x$ is not a marked point then $n_x\equiv 0 \pmod{p}$.
\item
    If $x$ is supersingular then $n_x=0$.
  \end{enumerate}
\end{prop}

\proof Suppose first that $x=x_i$. If $\nabla$ has logarithmic monodromy at
$x_i$, then $n_{x_i}=0$ by the discussion preceeding Proposition
\ref{nilpotprop1}. This proves (i). Suppose that $\nabla$ has toric monodromy
at $x_i$.  It is well known that the order of vanishing of a horizontal
section of a flat line bundle in a regular singular point is congruent to
minus the local exponent in that point. Since $\alpha_i$ is the local exponent
of $\M$ at $x_i$, we have that
\[
          n_{x_i} \;\equiv \; -\alpha_i \pmod{p}.
\]
The claim $n_{x_i}\not\equiv 0$ follows now from Proposition
\ref{indiprop1} (ii). This proves (ii). The proof of (iii) is similar.

Let $x$ be a supersingular point, and let $t$ be a local parameter at $x$. Put
$D=\partial/\partial t$. We have seen that we may regard $\Psi_\E(D^{\otimes
  p})$ as homomorphism of $\OO_X=\LL\to\M$. Let $e_0$ be the canonical
section. Then $\gamma(e_0)$ is a section of $\OO_X$ which is invertible at $x$.
It follows from (\ref{FLeq2}) that $\Psi_\E(D^{\otimes p})e_0=-D^{\otimes
  p}$. This shows that $n_x=0$. 
\Endproof

\section{Relation with deformation data}  \label{defodat}

\subsection{}\label{defodat1}

We start this section by recalling the notion of a deformation
datum. These objects arise in the study of bad reduction of Galois
covers of curves, see \S \ref{stablesec} for a quick review. We refer
to \cite{special} and \cite{bad} for a more thorough discussion.

\begin{defn} \label{defodatdef}
  A {\em deformation datum} on $X$ is a pair $(Z,\omega)$, where $Z\to
  X$ is a finite, at most tamely ramified Galois cover $\pi:Z\to X$ of
  smooth connected curves, and $\omega$ is a
  rational section of $\Omega_{Z/k}$. In addition, we require that the
  following holds.
  \begin{enumerate}
  \item
    Let $H$ be the Galois group of $Z\to X$.  For each $\sigma\in
    H$ we have $\sigma^*\omega=\chi(\sigma)\cdot\omega$, where
    $\chi:H\inj\FF_p^*$ is an injective character.
  \item
    The differential $\omega$ is logarithmic, i.e.\ of the form $\diff
    u/u$.
  \end{enumerate}
  Two deformation data $(Z,\omega)$ and $(Z',\omega')$ are called {\em
  equivalent} if there exists an isomorphism of $X$-schemes
  $\varphi:Z\iso Z'$ and an element $\epsilon\in\FF_p^\times$ such
  that $\varphi^*\omega'=\epsilon\cdot\omega$.  For each closed point
  $x\in X$ we define the following invariants:
  \begin{equation}\label{signatureeq}
       m_x := |H_z|, \qquad h_x := \ord_z\omega+1, \qquad 
       \sigma_x := h_x/m_x.
  \end{equation}
  Here $z\in Z$ is any point above $x$ and $H_z\subset H$ denotes the
  stabilizer of $z$. 
\end{defn}

\begin{defn}\label{ddssdef}
Let $(Z,\omega)$ be a deformation datum on $X$. 
\begin{enumerate}
\item A point $x\in X$
is said to be a {\em supersingular point} for $(Z,\omega)$ if
$\sigma_x=(p+1)/(p-1)$. 
\item A point $x\in X$ is said to be {\em singular}
if it is {\em not} a supersingular point and
$\sigma_x\not\equiv 1\pmod{p}$. 
\end{enumerate}
\end{defn}

The condition of Definition \ref{ddssdef} (ii) makes sense by Lemma
\ref{defodatlem} (i) below. Lemma \ref{defodatlem} (iii) implies that
there are only finitely many supersingular and singular points. We
refer to Example \ref{Legendreexa} for a motivation of the terminology
`supersingular'.

\begin{lem} \label{defodatlem}
  Let $(Z,\omega)$ be a deformation datum. 
  \begin{enumerate}
  \item For all $x\in X$, $m_x$ divides $p-1$. Moreover,
    $h_x$ and $m_x$ are relatively prime.
  \item
    If $h_x\not=0$ then $\gcd(p,h_x)=1$.
  \item
    For all but finitely many points $x\in X$ we have $\sigma_x=1$.
  \item
    We have
    \[
          \sum_{x\in X} (\sigma_x-1) = 2g-2.
    \]
  \end{enumerate}
\end{lem}

\proof The statement that $m_x$ divides $p-1$ follows immediately from
Definition \ref{defodatdef} (i). The statement that $\gcd(h_x, m_x)=1$ follows
from the assumption that $\chi$ is injective. 
Compare to the proof of \cite{special}, Prop.\ 2.5. Part (ii) is proved in
\cite{Henrio}, Cor.\ 1.8.b. Note that $h_x$ is denoted by $-m$ in that
paper. Part (iii) is obvious. Part (iv) follows from the Riemann--Roch Theorem.
\Endproof

\begin{defn} \label{signaturedef}
  An {\em signature} is given by a finite set $M$ and a map
  \[
      \bsigma: M \to \frac{1}{p-1}\cdot\ZZ,\quad x \mapsto \sigma_x
  \]
  such that $\sigma_x\geq 0$ and $\sigma_x\neq 1,\frac{p+1}{p-1}$ for
  all $x\in M$ and such that the number
  \[
        d := \frac{p-1}{2}\Big(\,2g-2
        -\sum_{x\in M} (\sigma_x-1) \,\Big)
  \]
  is a positive integer. The {\em singularities} of $\bsigma$ are the
  elements $x\in M$ with $\sigma_x\not\equiv 1\pmod{p}$. 

  Given a deformation datum $(Z,\omega)$ on $X$, the invariants
  $\sigma_x$ defined in \eqref{signatureeq} give rise to a signature
  $\bsigma$ (where $M$ is the set of points $x\in X$ with
  $\sigma_x\neq 1,\frac{p+1}{p-1}$). It follows from Lemma
  \ref{defodatlem} (iv) that the number $d$ defined above is the number of
  supersingular points.
\end{defn}

We shall always enumerate the singularities of a deformation datum
$(Z,\omega)$ as $x_1,\ldots,x_r$ and write $\sigma_i:=\sigma_{x_i}$
for $i=1,\ldots,r$.  We say that the deformation datum $(Z,\omega)$ is
{\em trivial} if $2g-2+r=0$.  In the rest of this paper, we exclude
trivial deformation data.

\subsection{}\label{stablesec}

We quickly recall how deformation data come up in the theory of
reduction of Galois covers. For more details, see e.g.\
\cite{bad}. Let $R$ be a complete discrete valuation ring with
fraction field $K$ of characteristic zero and residue field
$k=\bar{k}$ of characteristic $p>0$. Let $f_K:Y_K\to X_K$ be a
$G$-Galois cover defined over $K$.  We suppose, for simplicity, that
$p$ strictly divides the order of $G$ and that the curve $X_K$ has
good reduction, i.e.\ extends to a smooth $R$-curve $X_R$, with
special fiber $X:=X_R\otimes k$. Let $v$ denote the discrete valuation
of the function field of $X_K$ corresponding to the generic point of
$X$. It extends the valuation on $K$ and has residue field $k(X)$. Let
$A$ denote the completion of the valuation ring of $v$. Choose a
valuation of the function field of $Y_K$ extending the valuation $v$
via the map $f_K$, and let $C/A$ denote the resulting extension of
valuation rings. After replacing the ground field $K$ by some finite
extension, we may suppose that the ramification index of $C/A$ is one.

Suppose that the extension $C/A$ is ramified and let $B/A$ be the
maximal unramified subextension. The residue field of $B$ is the
function field $k(Z')$ of a Galois cover $Z'\to X$ of order prime to
$p$. It is shown in \cite{Raynaud98} that the extension $C/B$ is a
torsor under a finite flat group scheme over $R$. The special fiber of
this group scheme is isomorphic to either $\mu_p$ or $\alpha_p$. In
particular, the extension of residue fields of $C/B$ is inseparable of
degree $p$ and has the structure of a torsor under $\mu_p$ or
$\alpha_p$. It is a classical fact that such a torsor gives rise to a
differential form $\omega'\in \Omega_{k(Z')/k}$. If the group scheme is
$\mu_p$ (which we shall assume from now on) then $\omega'$ is
logarithmic. Namely, in this case the extension $C/B$ is defined by a
Kummer equation $y^p=h$, and the image of $h\in B$ in the residue
field $k(Z')$ is not a $p$th power. Then $\omega':=\diff h/h$. One
easily checks that the Galois group of $Z'/X$ acts on $\omega$ via a
character $\chi:H\to\FF_p^\times$. We see that the pair $(Z',\omega')$
has all the properties required for a deformation datum, except that
$\chi$ may not by injective. However, $\omega'$ descends to a
logarithmic differential form $\omega$ on the quotient $Z:=Z'/{\rm
Ker}(\chi)$, and we obtain a deformation datum $(Z,\omega)$.

The deformation datum $(Z,\omega)$ contains much valuable information
on the stable reduction of $Y_K$. The situation is particularly nice
if the cover $f_K:Y_K\to X_K$ is a {\em Belyi map}, i.e.\ if
$X_K=\PP^1_K$ and $f_K$ is branched precisely over the three points
$0,1,\infty$. If $f_K$ has bad reduction then the above construction
does always yield a deformation datum $(Z,\omega)$ (i.e.\ the
$\alpha_p$-case does not occur). Furthermore, the deformation datum
$(Z,\omega)$ is {\em special} (in the sense of Definition
\ref{specialdef} below) and almost completely determines the structure
of the stable reduction of the curve $Y_K$. See \cite{bad} for
details.

\begin{exa}\label{Legendreexa}
The Belyi map
\[
       \lambda:X(2p) \,\to\, X(2)\simeq\PP^1_\lambda
\]
between the modular curves of level $2p$ and $2$ is a Galois cover
branched at three points, with Galois group $G=\PSL_2(p)$. This cover
has a natural model over the field $\QQ(\zeta_p)$ which has bad
reduction to characteristic $p$.  As explained above, the reduction of
$X(2p)\to X(2)$ gives rise to a deformation datum $(Z, \omega)$. 

Let $p$ be an odd prime and $k$ an algebraic closure of the prime
field $\FF_p$. The {\em Hasse polynomial} is defined as
\[
     u := \sum_{i=1}^{(p-1)/{2}} 
         \binom{(p-1)/{2}}{i}^2 \lambda^i 
            \in k[\lambda].
\]
It has the property that the elliptic curve
\[
    E_\lambda:\qquad y^2=x(x-1)(x-\lambda)
\]
defined over $k$ is supersingular if and only if $u(\lambda)=0$.  It
is shown in \cite{mcav} that the deformation datum $(Z, \omega)$
associated to the cover $X(2p)\to X(2)$ is given by
\begin{equation}\label{Legendreeq}
      \omega := \frac{z\,\diff\lambda}{\lambda(\lambda-1)},
\end{equation}
and where $Z$ is the cover of $X=\PP^1_\lambda$ defined generically by
the equation $z^{(p-1)/2}={u}$. The crucial fact used in the proof is
that the Hasse polynomial $u$ is a solution to the Gauss
hypergeometric differential equation
\[
      \lambda(\lambda-1)\,u''+(2\lambda-1)\,u'+\frac{1}{4}\,u = 0.
\]
\end{exa}

In the rest of this section we give a correspondence between
deformation data and indigenous bundles in characteristic $p$. In
Example \ref{Legendreexa} above, the indigenous bundle corresponding
to the deformation datum $(Z,\omega)$ is essentially equivalent to the
Gauss hypergeometric differential equation. The construction we shall
give may seem very ad hoc. To fully appreciate it, one should have a
look at Mochizuki's theory of $p$-adic uniformization of ordinary
hyperbolic curves \cite{Mochizuki1}, \cite{Mochizuki2}, and the work
of Ihara on congruence relations \cite{Ihara74}.

\subsection{}  \label{defodat2}

We choose a rational function $t$ on $X$ with ${\rm d}t\neq 0$, and
set $D:=\partial/\partial t$.

\begin{lem}\label{defodatlem3}
  Let $Z\to X$ be a cyclic cover given by an
  equation $z^{p-1}=v^{-1}$, where $v\in k(X)$. Let
\[
\omega=z\, {\rm d}t.
\]
\begin{enumerate}
\item The differential form $\omega$ on $Z$ is logarithmic if and only
  if
\[
D^{p-1}(v)=-1.
\]
\item 
Suppose that $(Z,\omega)$ is a deformation datum, and let $t$
be a local parameter at $x\in X$. Then
\[
\ord_x(v)=(p-1)(1-\sigma_x).
\]
\end{enumerate}
\end{lem}

\proof It is a classical fact that a differential form $\omega=z\,
{\rm d}t$ is logarithmic if and only if
\[
D^{p-1} z=-z^p.
\]
For an outline of the proof see \cite{GS}, Exercise 9.6. This exercise is
stated only in the case that the genus of $X$ is zero, but one may easily
extend  the proof to our situation by considering $z$ as in element
of $k(\!(t)\!)$.
We deduce that
\[
D^{p-1} v= D^{p-1}\frac{z}{z^p}=\frac{1}{z^p}D^{p-1} z=-1.
\]
This proves (i). Part (ii) follows immediately from the definition of the
signature (\ref{signatureeq}).
\Endproof

For the rest of this subsection, we fix a deformation datum
$(Z,\omega)$ on $X$. It is clear that we may write $\omega=z\,\diff
t$, where $z$ is a rational function on $Z$ satisfying an equation
$z^{p-1}=v^{-1}$ with $v\in k(X)$. Set 
\[
\phi_0:=\omega^{\otimes(1-p)}=v\cdot D^{\otimes(p-1)}.
\]
Note that $\phi_0$ depends only on and determines the equivalence
class of the deformation datum $(Z,\omega)$. For each point $x\in X$
we define
\begin{equation} \label{nxeq}
   n_x :=\; \begin{cases}
      \;\; 0                \;&\; \text{\rm if $x$ is supersingular,}\\
      \;\;(p-1)\sigma_x     \;&\; \text{\rm if $x=x_i$ is singular,} \\
      \;\;(p-1)(\sigma_x-1) \;&\; \text{\rm otherwise}.
                \end{cases}     
\end{equation}
Note that $n_x$ is a nonnegative integer and is zero for all but
finitely many points $x$. Moreover, $n_x\equiv 0\pmod{p}$ if $x$ is
not singular.  Put $S:=\sum_xn_x\cdot x$. 

\begin{lem}\label{defodatlem4}
  The tensor $\phi_0$, considered as a rational section of
  the line bundle $(\olog_{X/k})^{\otimes(1-p)}(S)$, has no zeros.
  The poles of $\phi_0$ are precisely the supersingular points. All poles
  of $\phi_0$ are double poles.
\end{lem}

\proof Lemma \ref{defodatlem3} (ii) together with the definition (\ref{nxeq})
of $n_x$ implies that 
\[
  \ord_x v=(p-1)(1-\sigma_x)=\left\{\begin{array}{ll}
  -2 &\mbox{if $x$ is supersingular},\\
  p-1-n_x&\mbox{if $x$ is singular},\\
  -n_x&\mbox{otherwise}.
  \end{array}\right.
\]
It follows immediately that $\phi_0$, considered as a rational section
of the line bundle $(\olog_{X/k})^{\otimes(1-p)}(S)$, does not have
zeros and has poles of order $2$ in the supersingular points and no
poles elsewhere.  \Endproof

We now  associate to $(Z,\omega)$,  an indigenous bundle
$(\E,\nabla)$ on $(X;x_i)$ which is active, nilpotent and normalized.
Let $U\subset X$ denote the complement of the set of supersingular
points.  On $U$, we define a flat vector bundle $\E|_{U}$ of rank
two as follows:
\[
         \E|_{U} := \T(S)|_{U} \oplus \OO_{U}\cdot e_0.
\]
Let $t$, $D=\partial/\partial t$ and $\phi_0=v\cdot D^{\otimes(p-1)}$
be as above. The connection $\nabla$ is defined by the matrix equation
\[
   \nabla(D)(\underline{e}) = 
     \begin{pmatrix} 0 & v \\ 0 & 0 \end{pmatrix}\cdot\underline{e}
\]
with respect to the generic basis $\underline{e}=(D^{\otimes p},e_0)$
of $\E|_{U}$. It is clear that the definition of
$(\E|_{U},\nabla)$ depends only on $\phi_0$ and not on
the choice of the function $t$. We claim that $\nabla$ has regular
singularities in the singular points and no singularities
elsewhere. Namely, let $x\in U$ and assume that $t$ is a local
parameter at $x$. Set $D_x:=\partial/\partial t$ if $x$ is not a
singularity of the deformation datum and $D_x:=t\partial/\partial t$
otherwise. Then 
\[
    \nabla(D_x)(\underline{e}) = A\cdot\underline{e},
\]
where the matrix $A$ has coefficients which are regular in $x$. It
follows from Lemma \ref{defodatlem3} (i) that the $p$-curvature is
given by the rule
\begin{equation}\label{ddpcurveq}
   \Psi_\E(D^{\otimes p})(e) =D^{p-1}(v)\cdot\gamma|_{U}(e)\cdot
   D^{\otimes p}= -\gamma|_{U}(e)\cdot D^{\otimes p}.
\end{equation}
Here $\gamma|_{U}:\E|_{U}\to\OO_{U}$ is the projection onto the second
factor (the coefficient of $e_0$) and $e$ is an arbitrary 
section of $\E|_{U}$.

\begin{prop} \label{defodatprop} 
\begin{enumerate}
\item
  The bundle $(\E|_{U},\nabla)$ extends to an indigenous bundle
  $(\E,\nabla)$ on $X$. This extension is unique if we require that
  the section $e_0$ lies in the Hodge filtration $\Fil\E$ and that the
  projection $\gamma|_{U}$ extends to an epimorphism
  $\gamma:\E\to\OO_X$.
\item The bundle $(\E, \nabla)$ is active, nilpotent, and normalized.
\item The supersingular points (Definition \ref{ssdef}) are exactly the
  supersingular points for $(Z, \omega)$ (Definition \ref{ddssdef}). The
  marked points are the singularities of $(Z, \omega)$ (Definition
  \ref{ddssdef}).
\item The order $n_x$ of a spike is defined by (\ref{nxeq}).
\end{enumerate}
\end{prop}

\proof First we extend the bundle $\E|_{U}$ to the supersingular
points in such a way that the connection $\nabla$ is regular in these
points. This is a purely local problem. Fix a supersingular point $x$ for $(Z,
\omega)$.
We may suppose that the function $t$ introduced above is a local
parameter for $x$. We have $D^{p-1}(v-t^{p-1})=0$. Therefore, there
exists a rational function $w$ on $X$ such that $D(w)=v-t^{p-1}$. Since $v$
has a double pole at $x$, we may assume that $w$ has a simple pole.
Set $e_1:=e_0-w\cdot D^{\otimes p}$. We extend the bundle $\E|_{U}$
to the point $x$ in such a way that the pair
$\underline{e}':=(D^{\otimes p},e_1)$ is an $\OO_{X,x}$-basis of the
stalk $\E_x$. In terms of this basis, the connection $\nabla$ is given
by the matrix equation
\[
  \nabla(D)(\underline{e}') = 
   \begin{pmatrix} 0 & t^{p-1} \\ 0 & 0 \end{pmatrix}\cdot\underline{e}'.
\]
We see that $\nabla$ is regular in $x$. It is also clear that $\gamma$
extends to an epimorphism in the point $x$. This finishes the
definition of the flat bundle $(\E,\nabla)$.

We now  show that $(\E,\nabla)$ is indigenous. We define the
subbundle $\Fil\E$ as the maximal subbundle of rank one which has
$e_0$ as rational section. On $U$, the Kodaira--Spencer map 
\[
     \kappa: \Fil\E|_{U}=\OO_{U}\cdot e_0 \To 
      (\E/\Fil\E)|_{U}\otimes\olog_{U/k}= 
        \T(S)|_U\otimes\olog_{U/k}
\]
sends $e_0$ to $\phi_0$. Since $\phi_0$ has no zeros on $U$ (Lemma
\ref{defodatlem4} (i)), $\kappa$ is an isomorphism on $U$. Let $x$ be a
supersingular point. Using the notation of the preceeding paragraph, the stalk
$\Fil\E_x$ is generated by $e_0':=w^{-1}e_0$. Also, the pair $(e_0',e_1)$ is a
basis for the stalk $\E_x$. A short computation gives:
\begin{equation} \label{defodateq5}
    \nabla(D)(e_0') = w^{-1}t^{p-1}\cdot e_0' \,-\, 
      w^{-2}v\cdot e_1.
\end{equation}
Since $v$ has a pole of order $2$ and $w$ a simple pole in $x$, the
coefficient of $e_1$ in \eqref{defodateq5} has order $0$ in $x$. This means
that the Kodaira--Spencer map does not vanish in $x$. We have shown that
$(\E,\nabla)$ is indigenous.  

It is also clear that the extension we have defined is unique, under the
conditions that it is indigenous with $e_0$ lying in the Hodge filtration and
such that $\gamma|_{U}$ extends to an epimorphism $\gamma:\E\to\OO_X$. This
concludes the proof of  (i). 

Part (ii) follows from the construction of $(\E, \nabla)$.  Part (iv)
follows from (\ref{ddpcurveq}). Let $x\in X$ be a supersingular point
for $(Z, \omega)$ (Definition \ref{ddssdef}). The map $\gamma|_{U}$
extends to a map $\gamma:\E\to \OO_X$, by sending $e_0$ to $1$. This
makes $e_0$ into a canonical section (with respect to $\gamma$), as
defined in \S \ref{defodat2}. By definition of the extension of $\E$
to $x$, the section $e_1$ is invertible at $x$. This implies that
$e_0$ has a simple pole at $x$. Therefore $x$ is a supersingular
point of the bundle $(\E, \nabla)$. Since $e_0$ is regular on $U$, it
follows that the two notions of supersingularity agree.

Equation (\ref{nxeq}), together with Definition \ref{ddssdef}, implies that
$n_x\equiv 0\bmod{p}$ if $x$ is not supersingular and not a singularity of
$(Z, \omega)$. This implies that the marked points of $(\E, \nabla)$ are
exactly the singularities of $(Z, \omega)$.  \Endproof

 \begin{rem}\label{sigmarem} Lemma  \ref{defodatlem} (ii) states that
   $\sigma_x=h_x/m_x$ with $m_x|(p-1)$ and $\gcd(h_x, p)=1$. Note that this is
   compatible with Proposition \ref{nilpotprop1}. Namely, if $x=x_i$ is a
   marked point of $(\E, \nabla)$ then $\sigma_i=n_i/(p-1)$ where $n_i$ is
   either $0$ or prime to $p$. If $x$ is not a marked point, then $n_x\equiv
   0\bmod{p}$ and $\sigma_x=n_x/(p-1)+1=(n_x+p-1)/(p-1)$ and $\gcd(n_x+p-1,
   p)=1$.
\end{rem}

\subsection{} \label{defodat3}

We  now  reverse the construction of \S \ref{defodat2}.
Let $(\E,\nabla)$ be an indigenous  bundle on $X$ with regular
singularities in the marked points $x_1,\ldots,x_r$. Suppose that
$(\E,\nabla)$ is active, nilpotent, and normalized (\S
\ref{indidefsec}). 

\begin{prop} \label{defodatprop2}
  There exists a deformation datum $(Z,\omega)$ with singular
  points $x_1,\ldots,x_r$ such that the indigenous bundle associated
  to $(Z,\omega)$ by the construction of \S \ref{defodat2} is
  isomorphic to $(\E,\nabla)$. Moreover, $(Z,\omega)$ is unique up to
  equivalence. 
\end{prop} 

\proof Choose a horizontal and surjective homomorphism
$\gamma:\E\to\OO_X$ and identify the kernel of $\gamma$ with the line
bundle $\T(S)$, where $S=\sum_x n_x\cdot x$ is the zero divisor of the
$p$-curvature. Let $e_0$ be the canonical section of $\E$, corresponding to
the choice of $\gamma$ (\S \ref{indidefsec}). Set
  \[
         \phi_0 := \nabla(e_0).
  \]
Since $\gamma(e_0)=1$, by definition, it follows that
$\gamma(\nabla(D)e_0)=\nabla(D)\gamma(e_0)=0$. Therefore $\phi_0$ is a
rational section of the bundle
$\T(S)\otimes\olog_{X/k}=(\olog_{X/k})^{\otimes(1-p)}(S)$. It is easy
to check that $\phi_0$ does not depend on the choice of $\gamma$. 

Let $t$ be a rational function on $X$ with $\diff t\neq 0$ and write
$\phi_0=v\,(\diff t)^{\otimes(1-p)}$. Let $Z$ be a connected
component of the nonsingular projective curve with generic equation
  \[
          z^{p-1} = v^{-1}.
  \]
Set $\omega:=z\, \diff t$. One checks that the equivalence class of the
pair $(Z,\omega)$ does not depend on the choice of the function $t$
(the notion of equivalence we use is that of Definition \ref{defodatdef}).

We claim that the pair $(Z,\omega)$ is a deformation datum.  Indeed,
by definition, the pair $(Z,\omega)$ satisfies (i) of Definition
\ref{defodatdef}. Therefore it suffices to show that $\omega$ is
logarithmic. Let $D:=\partial/\partial t$. Then $D^p=0$. By the
definition of $\phi_0$, we have $\nabla(D)(e_0)=v\cdot D^{\otimes
p}$. It follows that
\begin{equation} \label{FLeq5}
   \Psi_\E(D^{\otimes p})(e_0) = \nabla(D)^{p-1}(v\cdot D^{\otimes p})
          = D^{p-1}(v)\cdot D^{\otimes p}.
\end{equation}
Recall that we have chosen $e_0$ such that $\gamma(e_0)=1$. Therefore,
\eqref{FLeq2} and \eqref{FLeq5} imply
\[
    D^{p-1}(v) = -1.
\]
Lemma \ref{defodatlem3} (i) implies that
$\omega$ is a logarithmic differential form.  This proves that $(Z,\omega)$ is
a deformation datum.

Let $(\E',\nabla')$ be the indigenous bundle associated to the deformation
datum $(Z,\omega)$, by the construction of \S \ref{defodat2}. It is clear that
the restriction of $\E$ and of $\E'$ to the complement $U$ of the
supersingular points are isomorphic. Therefore, it follows from the uniqueness
part of Proposition \ref{defodatprop} (i) that $\E$ and $\E'$ are isomorphic.
This finishes the proof of the lemma.  \Endproof

We can summarize the results of this section as follows.

\begin{thm} \label{defodatthm}
  The construction of \S \ref{defodat2} establishes a bijection
  between the set of equivalence classes of indigenous bundles on $X$
  which are active and nilpotent and the set of equivalence classes of
  deformation data. Suppose that the deformation datum $(Z,\omega)$
  corresponds to the indigenous bundle $(\E,\nabla)$. Then 
\begin{enumerate}
\item the supersingular points of $(Z, \omega)$ agree with the supersingular
  points of $(\E, \nabla)$,
\item the singular points of $(Z, \omega)$ are the marked points of $(\E,
  \nabla)$,
\item the points $x$ with $\sigma_x\neq 1$ and $\sigma_x\equiv 1\pmod{p}$ are
  the unmarked spikes of $(\E, \nabla)$ 
\end{enumerate}
The relation between the order $n_x$ of a spike $x$ of $(\E, \nabla)$ and the
invariant $\sigma_x$ associated to $(Z, \omega)$ is expressed by
(\ref{nxeq}).
\end{thm}

In the situation of Theorem \ref{defodatthm}, we may talk about the
signature of the deformation datum $(Z,\omega)$ as the signature of
the indigenous bundle $(\E,\nabla)$.

\begin{defn}\label{admdef}
Let $(X; \{x_1, \ldots, x_r\})$ be a marked curve over an
algebraically closed field $k$ of characteristic $p>0$. An indigenous
bundle $(\E, \nabla)$ on $(X; x_i)$ is called {\sl admissible} if it
does not have spikes outside the marked points.
\end{defn}

\begin{cor}\label{admcor}
  Let $(\E, \nabla)$ be an active and nilpotent indigenous bundle on $X$.  Let
  $r$ the number of marked points and $d$ the number of supersingular
  points. Suppose that 
  \[
     \frac{2d}{p-1}+\sum_{i=1}^r\sigma_i=2g(X)-2+r.
  \]
  Then $(\E, \nabla)$ is admissible. 
\end{cor}

\proof This follows immediately from Theorem \ref{defodatthm} and
Lemma \ref{defodatlem} (iv). In fact, the assumption implies that
there are no points $x\in X$ with $\sigma_x\neq 1$ except for the
singular and the supersingular points.  \Endproof


\section{Explicit construction of deformation data}
\label{desec}

In this section we suppose that $X=\PP^1_k$. We are interested in
constructing {\em deformation data} via their associated indigenous
bundles. It turns out that on $\PP^1$ it is more practical to replace
indigenous bundles by their associated differential equation.

\subsection{}\label{smssec}

Let $t$ denote the standard parameter on $X=\PP^1_k$. Set
$D:=\partial/\partial t$. We  write $f'$ instead of $D(f)$. Let $r\geq
3$ and $x_1,\ldots,x_r\in X$ pairwise distinct closed points. We
assume that $x_r=\infty$. Set $U_0:=X-\{\infty\}=\AA^1_k$ and
$U:=X-\{x_1,\ldots,x_r\}$.

Let $(\E, \nabla)$ be an indigenous bundle on $(X;x_i)$ which is active and
nilpotent. We assume that $(\E,\nabla)$ is normalized and choose a horizontal
and surjective morphism $\gamma:\E\to\OO_X$ (\S \ref{indidefsec}). We let
$\alpha_i\in\FF_p$ denote the local exponent at $x_i$ of the flat subbundle
${\rm Ker}(\gamma)\subset\E$. Then the local exponents at $x_i$ of
$(\E,\nabla)$ are $(\alpha_i,0)$. Denote by $d\geq 0$ the number of
supersingular points of $(\E,\nabla)$ (\S \ref{pcurvsec}).

The Hodge filtration $\Fil\E\subset\E$ is a line subbundle of $\E$. Let $e_1$
be an everywhere invertible section of $\Fil\E$ on $U_0$; it is unique up to
multiplication by a constant in $k^\times$.  Set $e_2:=\nabla(D)(e_1)$. It
follows from Definition \ref{indidef} (i) that ${\mathbf e}:=(e_1, e_2)$ forms
a basis of $\E$ on $U$. Hence we may write $\nabla(D)(e_2)=-p_2e_1-p_1e_2$, or
\begin{equation} \label{nablaeq}
    \nabla(D)(\mathbf{e}) = 
      \begin{pmatrix} 0 & -p_2 \\ 1 & -p_1 \end{pmatrix}\cdot\mathbf{e},
\end{equation}
where $p_1,p_2$ are rational functions which are regular on $U$. The
{\em differential operator} associated to $(\E,\nabla)$ is the second
order operator
\[
     L(u) := u''+p_1\,u'+p_2\,u.
\]

In \cite{Honda}, Appendix, one can find a more general dictionary
between flat bundles on $\PP^1$ and ordinary differential
operators. In particular, it is shown that the following properties of
the operator $L$ follow from the corresponding properties of
$(\E,\nabla)$:
\begin{itemize}
\item
  $L$ has regular singularities in the marked points $x_i$,
\item
  the local exponents at $x_i\neq \infty$  are $(\alpha_i,0)$ (the same as
  those of $(\E,\nabla)$).
\end{itemize}

\begin{rem}
  For some authors $L(u)=0$ is the differential equation associated to the
  {\em dual} of the bundle $(\E,\nabla)$, since it describes the horizontal
  vectors of the flat vector bundle dual to $(\E, \nabla)$, see \cite{Honda}
  appendix. The advantage of our convention is that the local exponents of
  $(\E, \nabla)$ coincide with the notion of local exponents of the
  differential operator $L$ in the classical sense, see for example
  \cite{Yoshida}, \S 1.2.5. 
\end{rem}

By \cite{Honda}, nilpotence of the $p$-curvature of $(\E,\nabla)$ is
equivalent to the assertion that the equation $L(u)=0$ has {\em
sufficiently many solutions in the weak sense}. This means that both
the equation $L(u)=0$ and its Wronskian equation $L_W(w)=w'+p_1w$ have
an algebraic solution. The nonvanishing of the $p$-curvature means
that the solution space of $L(u)=0$ has dimension one over $k(t^p)$.

\begin{prop} \label{Lprop1}
Let $L$ be the differential operator associated to the active, nilpotent,
normalized bundle $(\E, \nabla)$.
  \begin{enumerate}
  \item
  The equation $L(u)=0$ has a polynomial solution $u\in k[t]$ with
  simple zeros exactly in the supersingular points of $(\E,
  \nabla)$. It is unique up to multiplication by a constant in
  $k^\times$. 
  \item
    The local exponents of $L$ at $x_r=\infty$ are $(-d+\alpha_r,-d)$,
    where $d$ is the number of supersingular points of $(\E, \nabla)$.
  \end{enumerate}
\end{prop}

\proof It follows from the choice of $e_1$ that $u:=\gamma(e_1)$ is a
polynomial that does not vanish at the marked points. A simple
computation shows that $L(u)=0$. This proves the first assertion of
(i). The remaining assertions of (i) are clear. See also the proof of
Lemma \ref{sslem}. 

Write $(\gamma_1, \gamma_2)$ for the local exponents of $L$ at
$\infty$. An easy computation relating the local exponents of $L$ with those
of $(\E, \nabla)$ shows that $\gamma_1-\gamma_2=\alpha_r$ (up to
renumbering). The Fuchs' Relation (see \cite{Yoshida}, \S 1.2.6) states that
\begin{equation}\label{Fuchseq}
r-2=\gamma_1+\gamma_2+\sum_{i\neq r}\alpha_i=2\gamma_2+\sum_{i=1}^r
\alpha_i\in \FF_p.
\end{equation}

Let $n_x$ denote the spike order of $(\E,\nabla)$ at a point $x\in X$
and $\sigma_x$ the invariant attached to $x$ and the deformation datum
corresponding to $(\E,\nabla)$. Recall that
\[
\frac{n_x}{p-1}=\left\{\begin{array}{ll}
\sigma_x&\mbox{if }x=x_i,\\
0&\mbox{if $x$ is supersingular},\\
\sigma_x-1&\mbox{otherwise}.\end{array}\right.
\]
We have seen in \S \ref{monosec2} that $n_x\equiv -\alpha_i\bmod{p}$
if $x=x_i$ and is congruent to $0$ otherwise. Moreover,
$\sigma_x-1=2/(p-1)$ if $x$ is supersingular. Therefore Lemma
\ref{defodatlem} (iv) implies that
\[
     r-2 = \sum_{i=1}^r \sigma_i+\sum_{x\neq x_i} (\sigma_x-1)\equiv 
       -2d + \sum_{i=1}^r \alpha_i \pmod{p}.
\]
Together with \eqref{Fuchseq}, this shows that $\gamma_2=-d$ and
$\gamma_1=-d+\alpha_r$.  \Endproof

\subsection{}\label{accsec}

The goal of this section is to give a converse to Proposition
\ref{Lprop1}.  We fix distinct points $x_1, \ldots, x_r=\infty$ in
$\PP^1_k$ and a signature ${\boldsymbol\sigma}$ with $r$
singularities. Let $d$ be the positive integer of Defintion
\ref{signaturedef} (iii). We shall give a necessary and sufficient condition on a
differential operator $L$ to be associated to an indigenous bundle
with singularities $x_i$ and signature $\boldsymbol\sigma$, as in \S
\ref{smssec}.

Let $0\leq \alpha_i< p$ be such that
$\sigma_i\equiv\alpha_i\bmod{p}$. Put
\[
    Q=\prod_{i=1}^{r-1} (t-x_i)^{1-\alpha_i}.
\]
A general differential operator $L$ of degree two with regular singularities
at marked points $x_1, \ldots, x_r$ and local exponents $(\alpha_i,0)$ at
$x_i\neq\infty$ and $(-d+\alpha_r,-d)$ at $x_r=\infty$ can be written as
$L=(\partial/\partial t)^2+p_1(\partial/\partial t)+p_2$, where
\begin{equation} \label{p1p2eq}
     p_1 = \sum_{i=1}^{r-1} \frac{-\alpha_i+1}{t-x_i}, \qquad
     p_2 = \frac{d(d-\alpha_r) t^{r-3}-\beta_{r-4} t^{r-4}-\cdots
  -\beta_0}{\prod_{i=1}^{r-1} (t-x_i)}
\end{equation}
for constants $\beta_0,\ldots,\beta_{r-4}$. The constants $\beta_j$ are called
the {\em accessary parameters}.

Note that $Q'/Q=p_1$ by (\ref{p1p2eq}).
Before stating the existence result, we prove an easy lemma.

\begin{lem}\label{existencelem1}
  Let $(\E, \nabla)$ be an active and nilpotent indigenous bundle of
  signature ${\boldsymbol \sigma}$. Let $e_1,e_2$ be as in \S \ref{smssec}
  and $u$ as in Proposition
  \ref{Lprop1}. The horizontal sections of $(\E, \nabla)$ are
  \[
        Qw^p(-u'e_1+ue_2),
  \]
  for an arbitrary element $w\in k(t)$.
\end{lem}

\proof The assumption that $(\E, \nabla)$ is active implies that the solution
space is $1$-dimensional over $k(t)^p$. Therefore it suffices to show that
$L(Q(-u'e_1+ue_2))=0$.  This follows by direct verification.  \Endproof

\begin{prop}\label{existenceprop}
  Let $L=(\partial/\partial t)^2+p_1(\partial/\partial t)+p_2$ be a second
  order differential operator with regular singularities in $x_i$ and local
  exponents $(\alpha_i, 0)$ (resp.\ $(-d+\alpha_r, -d)$ in $x_i$ for $i\neq
  r$ (resp.\ $x_r$). Suppose that $L$ has a polynomial solution $u$ of degree
  $d$ such that
\begin{enumerate}
\item $u$ does not have zeros at the marked points,
\item we have
\[
   \ord_{x}\frac{1}{Qu^2}=\left\{\begin{array}{cl}
   0                     &\mbox{if $x$ is supersingular,}\\
   (p-1)\sigma_i-\alpha_i&\mbox{if }x=x_i,\\
   (p-1)(\sigma_x-1)&\mbox{otherwise}.
        \end{array}\right.
\]
\end{enumerate}
Then $L$ is associated to an indigenous bundle $(\E, \nabla)$ which is
active, nilpotent and normalized. Moreover, the signature of
$(\E,\nabla)$ is ${\boldsymbol \sigma}$. The supersingular points of
$(\E, \nabla)$ are exactly the zeros of $u$.
\end{prop}

\proof We define an indigenous bundle $(\E, \nabla)$ corresponding to $u$ as
follows.
We define a function $w$ by
\[
     D^{p-1}(\frac{1}{Qu^2})=-w^p.
\]
Set $\tilde{Q}:=w^{-p}Q$. Then
\begin{equation} \label{Qtildeeq}
     D^{p-1}(\frac{1}{\tilde{Q}u^2})=-1.
\end{equation}
On $U_0=\AA^1$ we let $\E$ be the trivial bundle with basis
$e_1,e_2$ and with the connection $\nabla$ defined by \eqref{nablaeq}.

A straightforward computation using the Fuchs' relation
(\ref{Fuchseq}) and expression for the local exponents at $\infty$
shows that we may take the pair $(e_1,\nabla(t\partial/\partial t)(e_1))$ as
a basis of $\E$ in a neighborhood of $x_r=\infty$, compare to
\cite{Yoshida}, \S 1.2.6. This defines the bundle $(\E,\nabla)$ on
$X=\PP^1_k$. We define the filtration $\Fil\E\subset\E$ as the maximal
line subbundle which contains $e_1$ as a rational section. It follows
immediately that the Kodaira--Spencer map is an isomorphism.

Let $\M$ be the kernel of the $p$-curvature of $\E$.  We leave it to
the reader to check that the image of $e_1/u$ in $\LL:=\E/\M$ is
invertible on $U_0$. One computes that
\[
\nabla(D) \frac{e_1}{u}=\frac{1}{\tilde{Q}u^2} \eta,
\]
where $\eta=\tilde{Q}(-u'e_1+ue_2)$ is a horizontal section of $\E$ (Lemma
\ref{existencelem1}). Therefore
\begin{equation} \label{Qteq}
\Psi_{\E}(D^{\otimes
  p})\frac{e_1}{u}=D^{p-1}(\frac{1}{\tilde{Q}u^2})\eta
     = - \eta,
\end{equation}
by \eqref{Qtildeeq}.  Using the notation of (\ref{FLeq2}), the
$p$-curvature considered as $\OO_X$-linear morphism
$\T\to\End_{\OO_X}(\E)$ is given by
\[
   \Psi_{\E}(D^{\otimes p})\frac{e_1}{u}=-D^{\otimes p}.
\]
Therefore, the bundle $(\E,\nabla)$ is indigenous, active, nilpotent
and normalized. Moreover, $e_0:=e_1/u$ is the canonical section of the bundle
$(\E,\nabla)$, as defined in \S \ref{pcurvsec}. 

It follows from \eqref{Qteq} that the order of zero of $\Psi_{\E}$ is
\[
n_x=\left\{\begin{array}{ll}
\alpha_i+\ord_xD^{p-1}(\frac{1}{Qu^2})&\mbox{ if $x=x_i$ is a marked point},\\
\ord_xD^{p-1}(\frac{1}{Qu^2})&\mbox{ otherwise}.
\end{array}\right.
\]
This shows that the signature of $(\E, \nabla)$ is as stated in the
proposition. 
\Endproof

Dwork's {\em accessary
parameter problem} is to find values for $(\beta_j)$ such that the
operator $L$ has nilpotent curvature. We are interested in those $L$
which satisfy the stronger conditions of Proposition \ref{existenceprop}
above. Namely, the $p$-curvature of $L$ should be nilpotent and
nonvanishing and $L(u)=0$ should have a polynomial solution of a given
degree $d$ which does not vanish at the marked points. 

\begin{defn}
  A {\em solution of the strong accessary parameter problem} with signature
  $(\sigma_i)$ and supersingular degree $d$ is a tuple $(x_i; \beta_j)\in
  k^r\times k^{r-3}$ such that the associated operator $L$ has a polynomial
  solution $u$ of degree $d$ which satisfies the conditions of Proposition
  \ref{existenceprop}.
\end{defn}

Let $(x_i; \beta_j)$ be a  solution of the strong accessary
parameter problem. Let $L$ be the differential operator and $(\E,\nabla)$ the
(normalized) indigenous bundle corresponding to the solution. Let $u$ be the
polynomial solution of $L(u)=0$ of degree $d$. According to Theorem
\ref{defodatthm}, $(\E,\nabla)$ gives rise to a deformation datum
$(Z,\omega)$. The following proposition shows that we can write down
$(Z,\omega)$ very explicitly in terms of the solution $u$. 

Let $x_{r+1},\ldots,x_s$ be the set of spikes which are not marked
points.  As usual, we write $n_i=n_{x_i}$ for the order of
the $p$-curvature $\Psi_\E$ at $x_i$.
Let $a_i$ be the unique integer such that $0\leq a_i<p-1$ and $a_i\equiv
n_i\bmod{p-1}$.  We may define a positive integer $\nu_i$ by
$\nu_i=(n_i-a_i)/(p-1)$ if $1\leq i\leq r$ and $\nu_i=1+(n_i-a_i)/(p-1)$
otherwise. Recall that $n_i\equiv -\alpha_i\bmod{p}$ if $1\leq i\leq r$ and
$n_i\equiv 0\bmod{p}$ otherwise. Therefore $\alpha_i\equiv
\nu_i-a_i\bmod{p}$ if $1\leq i\leq r$.

\begin{prop} \label{Zomegaprop}
  There is a constant $\epsilon\in k^\times$ such that the deformation datum
 is  $(Z,\omega)$,  where $Z\to X$ is the tamely
  ramified cyclic cover with generic equation
  \begin{equation} \label{Zeq}
      z^{p-1} = \prod_{i\neq r}(t-x_i)^{a_i}\,u^2
  \end{equation}
  and 
  \begin{equation} \label{omegaeq}
      \omega = \frac{\epsilon z\diff t}{\prod_{i\neq
      r}(t-x_i)^{-\nu_i+1}}.
  \end{equation}
\end{prop}

\proof At every point $x\in X$ we have the invariant $\sigma_x=h_x/m_x$
attached to the deformation datum $(Z,\omega)$. We have seen in \S
\ref{defodat} that $\sigma_x=(p+1)/(p-1)$ if and only if $x$ is a
supersingular point. By Proposition \ref{Lprop1}, this means that $x$ is a
zero of the polynomial $u$. At all other points $x=x_i$ we have
\[
   \sigma_i:=\sigma_{x_i} = \nu_i+\frac{a_i}{p-1}
\]
with $a_i$ and $\nu_i$ as defined before the statement of the
proposition.  It is now easy to see that the cyclic cover $Z\to X$ can
be identified with the cover defined by \eqref{Zeq}. Furthermore, the
differential in \eqref{omegaeq} has the same divisor as $\omega$ and
is therefore equal to $\omega$ for a suitable constant $\epsilon$. The
proposition is proved.  \Endproof

In \cite{IreneHabil}, Proposition 3.2.2 a direct proof of this proposition is
given which does not use the comparison with indigenous bundles. A key point
is to show that $1/Qu^2$ does not have residues in the supersingular
points. This follows in our set-up from Proposition
\ref{nilpotprop1} (iv). 

Given a deformation datum $(Z,\omega)$, its signature
${\boldsymbol \sigma}$ determines the number $d$ of
supersingular points, by Lemma \ref{defodat} (iv):
\begin{equation} \label{deq}
  d = \frac{p-1}{2}\big(s-2-\sum_{i=1}^s\sigma_i\big).
\end{equation}
The proof of Proposition \ref{Zomegaprop} shows that $(Z,\omega)$ is
essentially determined by its signature and the position of the
singular and the supersingular points. Therefore, in order to
construct a deformation datum with given signature, one
has to determine the position of the singular points $x_i$ and the
polynomial $u$ of degree $d$ (whose zeros are the supersingular
points) such that the differential $\omega'$ defined in Proposition
\ref{Zomegaprop} is logarithmic. Modulo the automorphisms of
$X=\PP^1$, this conditions gives rise to $2+d$ equations in
$r+d-3$ variables.

On the other hand, by Proposition \ref{Zomegaprop} the deformation
datum $(Z,\omega)$ corresponds to a solution of the strong accessory
parameter problem. This problem corresponds to $2(r-3)$ equations in
$2(r-3)$ variables. So if $d$ is large compared to $r$, this gives a
much better method to construct deformation data. This is particularly
striking for $r=3$, as the next example illustrates.

\begin{exa} \label{hgexa}
  Suppose that $r=3$. We may assume that the marked points are
  $(x_1,x_2,x_3)=(0,1,\infty)$.  Corollary \ref{admcor} implies that there are
  no spikes outside the marked points. 
 Equation \ref{deq} implies  that
  $\nu_i=0$ for $i=1,2,3$ and that $a_1+a_2+a_3\leq p-1$ is even. We
  have $d=(p-1-a_1-a_2-a_3)/2$. Since the number of accessary
  parameters is zero, there exists a unique order two differential equation
  with regular singularities at $0,1,\infty$ with local exponents
  $(-a_1,0),(-a_2,0),(-a_3-d,-d)$. This is the {\em hypergeometric
  equation} 
  \[
      t(t-1)\,u''+[(A+B+1)t-C]\,u'+AB\,u = 0,
  \]
  where $A:=(1+a_1+a_2+a_3)/2$, $B:=(1+a_1+a_2-a_3)/2$ and $C:=1+a_1$.
  It is proved in \cite{mcav}, Proposition 3.2, that this equation has
  a polynomial solution $u$ of degree $d$. It follows that \eqref{Zeq}
  and \eqref{omegaeq} define a deformation datum $(Z,\omega)$ with
  singular points $0,1,\infty$ and signature
  $(a_i/(p-1))_{i=1,2,3}$. 

  For $a_1=a_2=a_3=0$ this is Example \ref{Legendreexa}, giving
  rise to the Hasse polynomial and Gauss' hypergeometric equation.
\end{exa}

\subsection{}\label{specialsec}
In the rest of this section, we consider the case of special deformation data.
As explained in \S \ref{stablesec}, these correspond to the reduction of
$G$-Galois covers of $\PP^1$ branched at $3$ points. Here $G$ is some group
whose order is strictly divisible by $p$.  Let $(Z,\omega)$ be a deformation
datum on $X=\PP^1$, with singular points signature
$(\sigma_1,\ldots,\sigma_s)$. As in \S \ref{accsec}, we write
$\sigma_i=\nu_i+a_i/(p-1)$ with $\nu_i\geq 0$ and $0\leq a_i<p-1$.

\begin{defn}\label{specialdef}
  The deformation datum $(Z, \omega)$ is called {\sl special} if
  $\nu_i=0$ for exactly three indices $i$ and $\nu_i=1$ for the other indices.
\end{defn}

Note that if $(Z, \omega)$ is special then $\sigma_x\not\equiv 1\bmod{p}$, by
Corollary \ref{admcor}. In other words, the corresponding indigenous bundle
$(\E, \nabla)$ is admissible (Definition \ref{admdef}).

We consider the case of special deformation data with four
singularities.  We may assume that the singular points are $x_4=0, x_1=1,
x_2=\lambda, x_3=\infty$, for a variable $\lambda$, and that
$\nu_1=\nu_2=\nu_3=0$ and $\nu_4=1$. Let $p\geq 7$ and $(a_1, a_2, a_3,
a_4)=(0,0,1,3)$.  Put $d=(p-5)/2$. 

We start by describing the equations $F,G\in \FF_p[\lambda, \beta]$ whose zero
locus corresponds to the set of deformation data with signature $(0,0,
1/(p-1), (p+2)/(p-1))$, by using the conditions of Proposition
\ref{existenceprop}. Set $u=\sum_{i\geq 0} u_i t^i$, with $u_0=1$.  Assuming
that $L(u)=0$, we find a recursion
\begin{equation}\label{receq}
\lambda A(i)u_{i+1}=(B(i)+\beta)u_{i}-C(i)u_{i-1},
\end{equation}
with
\begin{equation*}\begin{split}
    A(i)&=(i+1)(i+a_4), \qquad B(i)=i^2(1+\lambda)+i(a_2+a_4+\lambda(a_1+a_4))
, \\
    C(i)&=(i-1-d)(i-1-d-a_3).
\end{split}\end{equation*}
 for the coefficients of $u$.
Since $p-a_4=p-3>d=(p-5)/2$, we have that
$A(i)\neq 0$, for $i=0, \ldots, d$. We let $F$ be the numerator of $u(d+1)$.
Since $C(d+1)=0$, it follows that $F=0$ implies that $L$ has a solution $u$ of
degree $d$. (In \cite{Beukers} we find an interpretation of this equation as
eigenvalue problem.) One easily sees that $\deg_\beta(F)=d+1$.

Let 
\[
w=\frac{1}{Qu^2}=\frac{1}{t^3(t-1)(t-\lambda)u^2}.
\]
One computes that
\[
{\rm Res}_{t=0}
w=u_0^2+(-2u_1u_0+3u_1^2-2u_2u_0+u_0^2)\lambda^2+(u_0^2-2u_0u_1)\lambda.
\]
The recursion (\ref{receq}) for the $u_i$ implies that
\[
u=1-\frac{\beta}{3\lambda}t+\frac{\beta^2-4(1+\lambda)\beta-6}{24\lambda^2}
t^2+\mbox{higher order terms}.
\]
We replace $u$ by $24\lambda^2 u$ and define $G$ to be the numerator of the
residue at $t=0$ of $w$, i.e.\
\[
G=
576\lambda^2+(480\beta+576)\lambda+864+60\beta^2+480\beta.
\]

The problem is now to show that for every $p\geq 7$ there exist solutions
$(\lambda, \beta)\in \bar{\FF}_p^2$ of $F=G=0$ such that the  points
$(0, 1, \infty, \lambda, \tau_i)$ are pairwise distinct. Here $(\tau_i)$ are
the zeros of $u$, i.e.\ the supersingular points. For $x\neq 0,1,\lambda,
\infty$, we have that $\ord_x(u)\equiv 1\bmod{p}$. Since $\deg_x(u)<p$, it
follows that if the points $(0,1,\lambda, \infty, \tau_i)$ are not
pairwise distinct, then $\lambda\in\{0,1,\infty\}$.

A solution $(\lambda, \beta)\in \bar{\FF}_p^2$ of $F=G=0$ is called a {\sl
  good solution} if $\lambda\not\in \{0,1\}$.
Table 1 gives a list of the degree of the total solution space and the degree
of the good solutions as function of $p$. We already remarked that all good
solutions have multiplicity one. The progression of the good degree strongly
suggests that there are deformation data for all $p\geq 7$.

\begin{table}[top]
\[
\renewcommand{\arraystretch}{1.5}
\begin{array}{c|c|c|c||c|c|c|c}
p&d&\mbox{total deg}&\mbox{good deg}
&p&d&\mbox{total deg}&\mbox{good deg}\\
\hline
7&1&12&3&53&29&148&48\\
11&3&24&8&59&27&168&56\\
13&4&30&10&61&28&174&58\\
17&6&23&3&67&31&191&62\\
19&7&48&16&71&33&203&67\\
23&9&60&18&73&24&150&50\\
29&12&78&25&79&37&226&73\\
31&13&83&27&83&39&239&79\\
41&18&112&36&83&39&239&79\\
43&19&118&38&97&46&282&93\\
47&21&132&43&&&&\\
\end{array}
\]
\caption{The number of deformation data with signature $(0,0,
    \frac{1}{p-1}, \frac{p+2}{p-1})$}
\end{table}

We finish the paper with another example. Let $p\geq 11$ and $(a_1, a_2, a_3,
a_4)=(0,0,p-8,3)$. We have that $d=(p-1-a_1-a_2-a_3-a_4)/2=2$.
We compute equations $F, G\in \FF_p[\lambda, \beta]$ for the locus of
deformation data, as above. One finds that
\begin{equation*}
\begin{split}
F&=
40\beta\lambda^2+(80\beta+14\beta^2+360)\lambda+132\beta+14\beta^2+360+\beta^3,
\\
G&= -4\beta\lambda^3+(3\beta^2-12\beta-24)\lambda^2+(12-8\beta)\lambda+12.
\end{split}
\end{equation*}
  
As in the previous example, the bad solutions satisfy
$\lambda=0,1$. Using $p\geq 11$, one shows that there are no trivial solutions
in this case. We conclude that for all $p\geq 11$ there exist special
deformation data with signature $(0,0,(p-8)/(p-1), (p+2)/(p-1))$.  

It is plausible but unproven that every special deformation datum
arises from a Belyi map with bad reduction. For results in this
direction see \cite{RRR}. In this paper one finds a sample of possible
applications of existence results for deformation data to Galois theory.

\section{Ordinary indigenous bundles} \label{ordinarysec}

\subsection{} 
In this section we discuss Mochizuki's notion of hyperbolically ordinary
indigenous bundles. For a discussion of this notion we refer to
\cite{Mochizuki2}, Introduction \S 2.1 (page 72) and Chapter IV.

Let $(g,r)$ be integers such that $2g-2+r\geq 3$.  Denote by $\M_{g,r}/\FF_p$
the stack of $r$-marked curves $(X; \{x_1, \ldots, x_r\})$ with $g(X)=g$. For
simplicity, we call such curves {\sl $(g,r)$-marked}. Let $\Nc_{g,r}$ be the
stack of admissible indigenous bundles $(\E, \nabla)$ on an $(g,r)$-marked
curve $(X; x_i)$. We write
\[
\pi:\Nc_{g,r}\To \M_{g,r}
\]
for the natural projection.

\begin{defn}\label{ordinarydef}
  Let $(\E, \nabla)$ be an admissible indigenous bundle on a $(g,r)$-marked
  curve $(X; x_i)$. We say that $(\E, \nabla)$ is {\sl hyperbolically
    ordinary} if $\pi:\Nc_{g,r}\to \M_{g,r}$ is \'etale at the point
  corresponding to $(\E, \nabla)$.  We write $\Nc_{g,r}^\ords\subset
  \Nc_{g,r}$ for the substack of hyperbolically ordinary indigenous bundles.
  A $(g,r)$-marked curve $(X; x_i)$ is called {\sl hyperbolically ordinary} if
  it lies in the image of $\pi:\Nc^\ords_{g,r}\to \M_{g,r}$.
\end{defn}

The notion hyperbolically ordinary should not be confused with the notion of an
ordinary abelian variety. Since the
latter notion plays no role in this paper, we will call ``hyperbolically
ordinary'' ordinary for short.

Example \ref{Legendreexa} implies that $(\PP^1_k; \{0,1,\infty\})$ is
ordinary. See also \cite{Mochizuki1}, page 1045 and \cite{Mochizuki2},
\S IV.2.1, page 205. Mochizuki shows that the generic $(g,r)$-marked
curve is ordinary, using the case $(g,r)=(0,3)$ and a deformation
argument (\cite{Mochizuki1}, Corollary III.3.8, page 1048). In this
section we want to show that
\[
\pi:\Nc^\ord_{0,4}\To \M_{0,4}
\]
is surjective. In other words, {\em every} $(0,4)$-marked curve is
ordinary. For this we use the translation of indigenous bundles in
terms of solutions of differential equations (\S \ref{accsec}). This
result could probably also be shown directly using the language of
indigenous bundles. Mochizuki shows that $\pi:\Nc_{g,r}\to \M_{g,r}$
is finite and flat of degree $p^{3g-3+r}$ (\cite{Mochizuki1}, Theorem
II.2.3, page 1029f). For $g=0$ this already follows from the work of
Dwork \cite{Dwork}.  For $(g,r)=(0,4)$ we find therefore that $\pi$ is
separable. This also follows from Lemma \ref{ordlem2} below.

\subsection{}\label{ordr=4sec}

Let $(\PP^1_k; x_i)$ be a $(0,4)$-marked curve. We may choose a
suitable coordinate $t$ on $\PP^1_k$ such that $x_1=0, x_2=1,
x_3=\lambda, x_4=\infty$.  This identifies $\M_{0,4}$ with
$\PP^1_\lambda-\{0,1,\infty\}$. Our goal is to construct indigenous
bundles $(\E, \nabla)$ on $(\PP^1_k; x_i)$ with logarithmic monodromy
in $x_i$ for $i=1, \ldots, 4$. Proposition \ref{existenceprop}
translates this into the existence of suitable solutions of the
differential operator
\begin{equation}\label{Lordeq}
L=\frac{\partial^2}{\partial
  t^2}+\left(\frac{1}{t}+\frac{1}{t-1}+\frac{1}{t-\lambda}\right)
\frac{\partial}{\partial t}+\frac{t-\beta}{t(t-1)(t-\lambda)},
\end{equation}
Equation (\ref{p1p2eq}).  Since the local exponents of $L$ at $t=\infty$ are
$(1, 1)$, it follows that if $L$ has a polynomial solution $u$ then the degree
of $u$ is $\equiv -1\bmod{p}$.

\begin{lem}\label{ordlem1}
\begin{enumerate}
\item Let $(\beta, \lambda)$ be such that $L=L_{\lambda, \beta}$ has a
  polynomial solution $u\in k(\lambda, \beta)[t]$. Then $L$ has a unique monic
polynomial solution of degree $p-1$.
\item For $L=L_{\lambda, \beta}$ as in (i), the corresponding indigenous
  bundle is  active and admissible (Definition \ref{admdef}).
\end{enumerate}
\end{lem}

\proof
Part (i) is proved in \cite{Beukers}, Lemma 1. 

Suppose that $L$ has a polynomial solution $u\in k(\lambda, \beta)[t]$ of
degree $p-1$.  Fix $i\in \{1,2,3\}$. Since the local exponents of $L$ at
$t=x_i$ are $(0,0)$, we have that $\ord_{x_i}(u)\equiv 0\bmod{p}$. Since the
degree of $u$ is $p-1<p$, we conclude therefore that $u$ does not have a zero
at $t=x_i$. Proposition \ref{existenceprop} implies therefore that $u$
corresponds to an indigenous bundle $(\E, \nabla)$.

Suppose that the $p$-curvature of $L$ is zero. Then the differential
operator $L$ has two linearly independent polynomial solutions $u_1,
u_2$ (\cite{Honda}, appendix.) Now Proposition 5.1 of \cite{Honda}
implies that the local exponents $(\gamma_1, \gamma_2)$ of $L$ at
$t=\infty$ are distinct. But this contradicts our assumptions. We
conclude that $\E$ is active. Corollary \ref{admcor} implies that
$(\E, \nabla)$ does not have unmarked spikes.  \Endproof

Beukers' proof of Lemma \ref{ordlem1} (i) is a concrete version of
\cite{Honda}, Lemma 1, page 174. The proof uses the recursion for the
coefficients of a solution $u$, together with a description of the accessary
parameter problem as eigenvalue problem. Beukers shows   that
$\pi:\Nc_{0,4}\to \M_{0,4}$ is surjective.  This also follows from the work of
Mochizuki. Unfortunately, Beukers' description only seems to works for $4$
marked points.

\begin{lem}\label{ordlem2}
  There exists an irreducible component
  $\Nc_1\subset \Nc_{0,4}$ such that the restriction $\pi_1:\Nc_1\to
  \M_{0,4}\simeq \PP^1_\lambda-\{0,1,\infty\}$ of $\pi$ has degree
  $1$.
\end{lem}
 
\proof Let $\Nc_1=\Spec(k[\lambda,1/\lambda(\lambda-1),
\beta])/(\lambda+\beta+1)$, and write $\pi_1:\Nc_1\to \PP^1_\lambda$ for the
natural projection. We claim that the differential operator $L$ has a unique
monic solution $u$ of degree $p-1$ over $\Nc_1$.

Let $\tilde{u}=\sum_{i\geq 0} u_i t^i\in \bar{\FF}_p(\Nc_1)[[t]]$ be a
solution of $L$ with $u_0\neq 0$. One easily checks that the
coefficients $u_i$ of $u$ satisfy a recursion
\begin{equation}\label{receq2}
\lambda(i+1)^2u_{i+1}=(1+\lambda)(i^2+i+1)u_i-i^2u_{i-1}.
\end{equation}
Note that $u:=\sum_{i=0}^{p-1} u_i t^i$ is a solution of $L$ if and only if 
\begin{equation}\label{acccondeq}
(1+\lambda)u_{p-1}-u_{p-2}=0.
\end{equation}
Moreover, $u_1, \ldots u_{p-1}$ are uniquely determined by (\ref{receq}) and
$u_0$. 

Put $v_{p-1-i}:=u_{i}\lambda^i$ for $i=0, \ldots, p-1$. One
easily deduces from (\ref{receq2}) that the $v_i$ satisfy the same recursion
(\ref{receq2}). Since the coefficients $u_0, \ldots u_{p-1}$ are uniquely
determined by $u_0$ and the recursion, it follows that 
\[
\frac{u_i}{u_0}=\frac{v_i}{v_0}=\frac{u_{p-1-i}}{\lambda^{p-1-i}u_{p-1}},
\qquad i=0, \ldots p-1.
\]
Since $u_1=(1+\lambda)u_0/\lambda$, we deduce that (\ref{acccondeq}) is
satisfied.
\Endproof

Lemma \ref{ordlem2} can be understood as follows. The marked curve
$(\PP^1_k; \{0,1,\lambda, \infty\})$ has an automorphism $\sigma:t\mapsto
\lambda/t$.  In Lemma \ref{ordlem2}, we have chosen the accessary parameter
$\beta$ so that $\sigma$ is a symmetry of the differential operator $L$.

\begin{prop}\label{ordprop}
The morphism
\[
\pi:\Nc^\ords_{0,4}\To \M_{0,4}
\]
is surjective.
\end{prop}

\proof Lemma \ref{ordlem2} implies that $\Nc_{0,4}$ has an irreducible
component such that the restriction of $\pi$ to this component has
degree $1$.  We denote this component also by $\Nc_1$. Lemma
\ref{ordlem1} (ii) states that all points of $\Nc_{0,4}$ corresponds
to admissible bundles. Therefore it follows from \cite{Mochizuki1},
Corollary II.2.16 (page 1043), that $\Nc_1$ is smooth over ${\FF}_p$.
Hence $\Nc_1$ does not intersect the rest of $\Nc_{0,4}$. This implies
that the restriction of $\pi$ to $\Nc_1$ is \'etale, and $\Nc_1$ is a
substack of $\Nc_{0,4}^\ords$. The proposition follows.  \Endproof

The following corollary immediately follows from the proof of Proposition
\ref{ordprop}. This answers an expectation of Mochizuki (\cite{Mochizuki2},
first remark on page 205).

\begin{cor}\label{ordcor}
The stack $\Nc_{0,4}$ is disconnected for all $p\geq 3$.
\end{cor}

It is interesting to study what the complement of $\Nc_1$ in $\Nc_{0,4}$ looks
like.  The following lemma describes this for $p=5$. It can be easily checked
using for example magma. The lemma illustrates that admissible and ordinary
are in general two different notions. 

\begin{lem}\label{ordlem3}
  Let $p=5$. Then $\Nc_{0,4}$ is the disjoint union of two smooth absolutely
  irreducible components $\Nc_1$ and $\Nc_2$. The restriction of $\pi$ to
  $\Nc_2$ has degree $4$, and  is ramified at 12 points of order $2$. 
The branch locus consists of $8$ points.
\end{lem}

\subsection{} \label{ordtorsec}

To finish this section, we consider the analogous question in the torally
indigenous case. Namely, for $i=1,2,3,4$ let $0\leq \sigma_i=a_i/(p-1)<1$ be
rational numbers, and put ${\boldsymbol \sigma}=(\sigma_i)$.  We write
$\Nc_{0,4}({\boldsymbol \sigma})$ for the stack of admissible indigenous
bundles of signature ${\boldsymbol \sigma}$. Here the signature of an
indigenous bundle is the signature of the corresponding deformation datum.
Recall that the signature can be expressed in terms of the order of the
spike in $t=x_i$. The stacks $\Nc_{0,4}({\boldsymbol \sigma})$ do not occur
in the work of Mochizuki, since Mochizuki does not consider the signature as
combinatorial  invariant. However, $\Nc_{0,4}({\boldsymbol \sigma})$ is a
substack of $\Nc_{g,r}[a_1+a_2+a_3+a_4]$ of admissible indigenous bundles
which are spiked of strength $a_1+a_2+a_3+a_4$ (\cite{Mochizuki2},
Introduction \S 1.2, page 39).

The following proposition is proved in \cite{IreneHabil},
\S 3.4. It uses the deformation theory of $\mu_p$-torsors (\cite{cotang}). The
proposition illustrates that the situation becomes more complicated for
arbitrary signature. Namely, the morphism
$\Nc_{0,4}^\ords({\boldsymbol\sigma})\to \M_{0,4}$ is not always surjective.

\begin{prop}\label{ordprop2}
 Suppose that $a_1+a_2+a_3+a_4$ is even and nonzero. We define a unique
  integer $0\leq d<p-1$ by the property
\[
2d+a_1+a_2+a_3+a_4\equiv 0\bmod{p},
\]
and put $\nu=(2d+a_1+a_2+a_3+a_4)/(p-1).$
\begin{enumerate}
\item If $\nu=1$ and $\Nc_{0,4}({\boldsymbol \sigma})$ is nonempty, then $\dim
  \Nc_{0,4}({\boldsymbol \sigma})=0$.
\item If $\nu=2$, then $\dim \Nc_{0,4}({\boldsymbol \sigma})=1$, and
  $\pi:\Nc_{0,4}^\ords({\boldsymbol \sigma})\to \M_{0,4}$ is surjective.
\item $\nu=3$, then $\Nc_{0,4}({\boldsymbol \sigma})$ is empty.
\end{enumerate}
\end{prop}

\proof The statement on the dimension of $\Nc_{0,4}({\boldsymbol \sigma})$ is
proved in \cite{IreneHabil}, Lemma 3.4.2. The surjectivity of $\pi$ in (ii) is
proved in \cite{IreneHabil}, Proposition 3.4.3. The nonexistence in case (iii)
follows from Lemma \ref{defodatlem} (iv).
\Endproof

A similar proposition for special deformation data has been proved in
\cite{cotang}, Theorem 5.14. Let ${\boldsymbol \sigma}=(\sigma_i)$ be
the signature of a special deformation datum. I.e.\ there are three
indices $i$ such that $0\leq \sigma_i<1$. For all other $i$ we have
that $1\leq  \sigma_i<2$. Then the space $\Nc_{0,
  3}({\boldsymbol \sigma})$ of special deformation data has dimension
$0$, if it is nonempty. If ${\boldsymbol \sigma}=(\sigma_1, \sigma_2,
\sigma_3)$ is as in Example \ref{hgexa}, then $\Nc_{0,
  3}({\boldsymbol \sigma})$ is exactly one point. Therefore, the
corresponding indigenous bundles could also be called ordinary.


\vspace{6ex}

\begin{minipage}[t]{5cm}
Mathematisches Institut\\
Heinrich-Heine-Universit\"at\\
Universit\"atsstr. 1\\
40225 D\"usseldorf \\
bouw@math.uni-duesseldorf.de\\
\end{minipage}
\hfill
\begin{minipage}[t]{5cm}
\begin{flushright}
Mathematisches Institut\\    
Universit\"at Bonn\\
Beringstr. 1\\
53115 Bonn\\
wewers@math.uni-bonn.de

\end{flushright}
\end{minipage}

\end{document}